\documentclass{article}\begin{document}
\centerline{\bf\large Formulas for 2-D Beltrami Operators}\vskip .2in

\centerline{Nikolaos D. Bagis}

\centerline{Aristotele University of Thessaloniki}
\centerline{Thessaloniki, Greece}
\centerline{nikosbagis@hotmail.gr}

\[
\]

\centerline{\bf Abstract}

\begin{quote}
In this article we consider 2-dimensional surfaces. We define some new operators which enable us to evaluate quantities of the surface, such invariants, in a more systematic way.      
\end{quote}


\section{Introduction}

I will use Pfaff differential forms. Assume a two-dimesional surface $S$ of the Eucledean space $E_3\cong \textbf{R}^3$ which is of class $C^3$. That is, the surface is given by 
$$
\overline{x}=\overline{x}(u,v)=\{x_1(u,v),x_2(u,v),x_3(u,v)\}\textrm{, }u,v\in D 
$$
and $x_{i}(u,v)\in C^3$,  $\overline{x}_u\times \overline{x}_v\neq \overline{0}$, $\overline{x}_u=\frac{\partial\overline{x}}{\partial u}$, $\overline{x}_v=\frac{\partial\overline{x}}{\partial v}$. 
In every point $P$ of the surface we attach a moving frame of three orthonormal vectors (that is $\{\overline{\epsilon}_1,\overline{\epsilon}_2,\overline{\epsilon}_3\}$ and  $\left\langle \overline{\epsilon}_i,\overline{\epsilon}_j\right\rangle=\delta_{ij}$), with the assumption that $\overline{n}=\overline{\epsilon}_3$ is orthonormal to the tangent plane of surface (in every $P$).Then there exist Pfaff (differentiatable) forms, $\omega_i$ and $\omega_{ij}$ such that
$$
d\overline{x}=\sum^{3}_{j=1}\omega_j\overline{\epsilon}_j\textrm{, }(\omega_3=0\Leftrightarrow \overline{n}=\overline{\epsilon}_3)
$$
$$
d\overline{\epsilon}_i=\sum^{3}_{j=1}\omega_{ij}\overline{\epsilon}_j\textrm{, }i=1,2,3
$$
This can be seen as:
$$
d\overline{x}=\{\partial_1 x_1du+\partial_2x_1dv,\partial_1 x_2du+\partial_2x_2dv,\partial_1 x_3du+\partial_2x_3dv\}
$$
and the Pfaff derivatives $\nabla_kf$ and $\nabla_k\overline{F}$ for any function $f$ or vector $\overline{F}$ (resp.) are defined as
$$
df=\sum^{3}_{k=1}(\nabla_kf)\omega_k=\sum^{3}_{k=1}\partial_kfdu_k.\eqno{(1)}
$$
Set now
$$
q_1=\frac{d\omega_1}{\omega_1\wedge\omega_2}\textrm{, }q_2=\frac{d\omega_2}{\omega_1\wedge\omega_2}.
$$
For to holds $(1)$ it must be
$$
\nabla_1\nabla_2f-\nabla_2\nabla_1f+q_1\nabla_1f+q_2\nabla_2f=0\textrm{, (condition)}.
$$
From the relations $d\left\langle\overline{e}_i,\overline{e}_j\right\rangle=0$, $d(d\overline{x})=\overline{0}$, $d(d\overline{n})=\overline{0}$, we get the structure equations of the surface:
$$
\omega_{ij}+\omega_{ji}=0\textrm{, }i,j=1,2,3,
$$
$$
d\omega_j=\sum^{3}_{i=1}\omega_i\wedge\omega_{ij}\textrm{, }j=1,2,3
$$ 
$$
d\omega_{ij}=\sum^{3}_{k=1}\omega_{ik}\wedge\omega_{kj}\textrm{, }i,j=1,2,3.
$$
Observe that $\omega_3=\omega_{11}=\omega_{22}=\omega_{33}=0$ and we can write
$$
d\overline{x}=\omega_1\overline{\epsilon}_1+\omega_2\overline{\epsilon}_2
$$
$$
d\overline{\epsilon}_1=\omega_{12}\overline{\epsilon}_2-\omega_{31}\overline{\epsilon}_3\eqno{(2)}
$$
$$
d\overline{\epsilon}_2=-\omega_{12}\overline{\epsilon}_1-\omega_{32}\overline{\epsilon}_3
$$
$$
d\overline{\epsilon}_3=\omega_{31}\overline{\epsilon}_1+\omega_{32}\overline{\epsilon}_2
$$
Moreover it is (structure equations):
$$
d\omega_1=\omega_{12}\wedge\omega_2
$$
$$
d\omega_2=-\omega_{12}\wedge\omega_1
$$
$$
\omega_1\wedge\omega_{31}+\omega_2\wedge\omega_{32}=0
$$
$$
d\omega_{12}=-\omega_{31}\wedge\omega_{32}\eqno{(3)}
$$
$$
d\omega_{31}=\omega_{12}\wedge\omega_{32}
$$
$$
d\omega_{32}=-\omega_{12}\wedge\omega_{31}
$$
If we write the connections (of $\omega_{ij}$ in terms of $\omega_i$):
$$
\omega_{12}=\xi\omega_1+\zeta \omega_2
$$
$$
\omega_{31}=-a\omega_1-b\omega_2
$$
$$
\omega_{32}=\eta\omega_1-c\omega_2
$$
We easily get (from the structure equations) $\xi=q_1$, $\zeta=q_2$, $\eta=-b$. Hence
$$
\omega_{12}=q_1\omega_1+q_2 \omega_2\eqno{(4)}
$$
$$
\omega_{31}=-a\omega_1-b\omega_2\eqno{(5)}
$$
$$
\omega_{32}=-b\omega_1-c\omega_2.\eqno{(6)}
$$
\\
\textbf{Definition 1.}\\
\textbf{i.} We call 2-Beltrami derivative the quantity
\begin{equation}
\Delta_2f:=\frac{d\left\langle\frac{df\wedge d\overline{x}}{\omega_1\wedge\omega_2},d\overline{x}\right\rangle}{\omega_1\wedge\omega_2}.
\end{equation}
\textbf{ii.} We call $f$ harmonic iff $\Delta_2f=0$.\\
\\

It also holds
\begin{equation}
\Delta_2f=\nabla_1\nabla_1f+\nabla_2\nabla_2f+q_2\nabla_1f-q_1\nabla_2f,
\end{equation}
where $\nabla_k$ are the directional Pfaff derivatives with respect to the diferental forms $\omega_1$, $\omega_2$ and
\begin{equation}
df=\left(\nabla_1f\right)\omega_1+\left(\nabla_2f\right)\omega_2
\end{equation}
Here the Pfaff differental forms $\omega_1$, $\omega_2$ are related with the surface $S:\overline{x}=\overline{x}(u,v)$ with the relation
\begin{equation}
d\overline{x}=\omega_1\overline{\epsilon}_1+\omega_2\overline{\epsilon}_2
\end{equation}
and $\omega_3=0$, (that is $\overline{\epsilon}_1$, $\overline{\epsilon}_2$ are in the tangent plane$-P$ of $S$ and $\overline{e}_3$ is vertical to $P$).\\

Note here that for every function $f$ holds the the following identity
\begin{equation}
\nabla_1\nabla_2f-\nabla_2\nabla_1f+q_1\nabla_1f+q_2\nabla_2f=0
\end{equation} 
\\
\textbf{Definition 2.}\\   
If for the functions $\Phi,\Phi^{*}:\textbf{R}\times\textbf{R}\rightarrow\textbf{R}$ hold the ralations
\begin{equation}
\nabla_1\Phi=-\nabla_2\Phi^{*}\textrm{ and }\nabla_{2}\Phi=\nabla_1\Phi^{*},
\end{equation}  
then we call $\Phi,\Phi^{*}$ analytic.

\section{Results}

\textbf{Theorem 1.}\\
If $\Phi$ is analytic, then $\Phi$ is also and harmonic.\\
\\ 
\textbf{Proof.}\\
From the definition of the analytic functions and identity (5) we get
$$
\Delta_2\Phi=\nabla_1\nabla_1\Phi+\nabla_2\nabla_2\Phi+q_2\nabla_1\Phi-q_1\nabla_2\Phi=
$$
$$
=-\nabla_1\nabla_2\Phi^{*}+\nabla_2\nabla_1\Phi^{*}-q_2\nabla_2\Phi^{*}-q_1\nabla_1\Phi^{*}=0
$$
\\
\textbf{Theorem 2.}\\
If $f$ is in $C^{(2)}$, then
\begin{equation}
\int_{\partial D}\frac{f}{\rho_g}ds+\int\int_{D}Kf\omega_1\wedge\omega_2=\int\int_{D}\left[q_2\nabla_1f-q_1\nabla_2f\right]\omega_1\wedge\omega_2+\int_{\partial D}fd\phi
\end{equation}
In case that $f$ is harmonic the above identity becomes
\begin{equation}
\int_{\partial D}\frac{f}{\rho_g}ds+\int\int_{D}\left(\nabla_1^2f+\nabla_2^2f+Kf\right)\omega_1\wedge\omega_2=\int_{\partial D}fd\phi
\end{equation}
\\
\textbf{Proof.}\\
It holds
\begin{equation}
\frac{ds}{\rho_g}=d\phi-\delta\psi
\end{equation}
where
\begin{equation}
\omega_{12}+\delta\psi=0\textrm{, }\omega_{12}=q_1\omega_1+q_2\omega_2
\end{equation}
and the Gauss curvature $K$ is
\begin{equation}
K=-\frac{d\left(q_1\omega_1+q_2\omega_2\right)}{\omega_1\wedge\omega_2}
\end{equation}
Hence
$$
\frac{f}{\rho_g}ds=fd\phi-f\delta\psi
$$
Integrating
\begin{equation}
\int_{\partial D}\frac{f}{\rho_g}ds=\int_{\partial D}fd\phi+\int_{\partial D}f\omega_{12}
\end{equation}
But from Stokes formula
$$
\int_{\partial D}f\omega_{12}=\int\int_{D}d[f\omega_{12}]=
$$
$$
=\int\int_{D}df\wedge \omega_{12}+\int\int_{D}fd\omega_{12}=
$$
$$
=\int\int_{D}\left[\nabla_1f\omega_1+\nabla_2f\omega_2\right]\wedge\left(q_1\omega_1+q_2\omega_2\right)-\int\int_{D}Kf\omega_1\wedge\omega_2=
$$
\begin{equation}
=\int\int_{D}\left[q_2\nabla_1f-q_1\nabla_2f\right]\omega_1\wedge\omega_2-\int\int_{D}Kf\omega_1\wedge\omega_2
\end{equation}
From (13) and (12) we get
\begin{equation}
\int_{\partial D}\frac{f}{\rho_g}ds=\int\int_{D}\left[q_2\nabla_1f-q_1\nabla_2f-Kf\right]\omega_1\wedge\omega_2+\int_{\partial D}fd\phi
\end{equation}
In case that $f$ is harmonic, using (2) and (5) we get the second result (relation (8)).\\
\\
\textbf{Theorem 3.}\\
Define
\begin{equation}
\Pi_2f:=q_2\nabla_1f-q_1\nabla_2f-Kf
\end{equation}
then
\begin{equation}
\int_{\partial D}\frac{f}{\rho_g}ds=\int\int_{D}\left(\Pi_2f\right)\omega_1\wedge\omega_2+\int_{\partial D}fd\phi
\end{equation}
\\
\textbf{Proof.}\\
Follows from the definition of $\Pi_2$ and (14).\\
\\
\textbf{Theorem 4.}\\
\begin{equation}
\Pi_2f=0\Leftrightarrow d\left(f\omega_{12}\right)=0 
\end{equation}
\\
\textbf{Proof.}\\
Assume that $\Pi_2f=0$, then for every $D$ such that $\partial D$ is simple connected closed curve, we have
$$
\int_{\partial D}f\left(\frac{ds}{\rho_g}-d\phi\right)=\int\int_{D}\left(\Pi_2f\right)\omega_1\wedge\omega_2=\int_{\partial D}f \omega_{12}=\int\int_{D}d\left(f \omega_{12}\right)=0
$$
\\
\textbf{Corollary 1.}\\
If the Gauss curvature of $S$ is $K=0$ then $S$ posess no closed geodesic curves.\\
\\
\textbf{Proof.}\\
Suppose that exist $K=0$ and $\int_{\partial D}\frac{ds}{\rho_g}=0$, then (7) with $f=1$ becomes
$$
\int_{\partial D}d \phi=0
$$  
which is no true since it must be $\int_{\partial D}d\phi=2\pi$.\\
\\ 
\textbf{Corollary 2.}\\
Let $A_{D}$ denote the area of a closed geodesic curve $\partial D$, then  
\begin{equation}
A_{D}=\int_{\partial D}\frac{d\phi}{K}+\int\int_{D}\left[q_2\nabla_1\left(\frac{1}{K}\right)-q_1\nabla_2\left(\frac{1}{K}\right)\right]\omega_1\wedge\omega_2
\end{equation}
when $K\neq0$ is the Gauss curvature of $S$.\\
\\
\textbf{Proof.}\\
Set $\Phi=\frac{1}{K}$ in (7).\\
\\
\textbf{Notes.}\\
\textbf{i.} If $S$ have constant non zero Gauss curvature and $\partial D$ is closed geodesic curve, then
\begin{equation}
A_{D}=\frac{2\pi}{K}=const
\end{equation} 
This result can be exctract from the well known integral formula of Gauss-Bonet:
\begin{equation}
\int_{\partial D}\frac{ds}{\rho_g}+\int\int_{D}K\omega_1\wedge\omega_2=2\pi
\end{equation} 
\textbf{ii.} In case we don't want to use $q_1$ and $q_2$ we can write (7) as  
\begin{equation}
\int_{\partial D}f\delta\psi +\int\int_{D}\left(\Delta_2f-\nabla_1^2f-\nabla_{2}^2f-Kf\right)\omega_1\wedge\omega_2=0
\end{equation} 
or the equivalent
\begin{equation}
\int_{\partial D}f\left(\frac{1}{\rho_g}-\frac{d\phi}{ds}\right)ds=\int\int_{D}\left(\Delta_2f-\nabla_1^2f-\nabla_{2}^{2}f-Kf\right)\omega_1\wedge\omega_2
\end{equation}
 
However in general holds the following relation
\begin{equation}
\omega_{12}=\left(\frac{1}{\rho_g}-\frac{d\phi}{ds}\right)ds=-\delta\psi,
\end{equation}
where $\delta\psi$ is the difference between the two angles, after a prallel transfer in the Levi-Civita sence in a closed path $\partial D$. It is known that 
\begin{equation}
\int_{\partial D}\delta\psi=\int\int_{D}K\omega_1\wedge\omega_2
\end{equation}
Hence
\begin{equation}
\int_{\partial D}f\delta\psi+\int\int_{D}\left(\Delta_2f-\nabla_1^2f-\nabla_{2}^{2}f-Kf\right)\omega_1\wedge\omega_2=0
\end{equation}
\\
\textbf{Proof.}\\
It holds 
\begin{equation}
\frac{1}{\rho_g}=\frac{\omega_{12}+d\phi}{ds}
\end{equation}
and
\begin{equation}
\left(\frac{1}{\rho_g}-\frac{d\phi}{ds}\right)ds=\omega_{12}=\widetilde{\omega}_{12}=-\delta\psi
\end{equation}
\\
\\
\textbf{Lemma 2.}\\
In every surface $S$ holds
\begin{equation}
\omega_{12}=\frac{ds}{\rho_g}-d\phi=\frac{ds^{III}}{\rho_g^{III}}-d\phi^{III}=\widetilde{\omega}_{12},
\end{equation}
where the symbol $III$ means the spherical image of $S$ i.e. the surface 
\begin{equation}
S^{III}:\overline{y}=\overline{\epsilon}_3=\overline{n}
\end{equation}
\\   
\textbf{Proof.}\\
One can easily see that
\begin{equation}
\widetilde{q}_1=\frac{d\omega_{31}}{\omega_{31}\wedge\omega_{32}}=\frac{q_2b-q_1c}{K}\textrm{ and }
\widetilde{q}_2=\frac{d\omega_{32}}{\omega_{31}\wedge\omega_{32}}=\frac{q_1b-q_2a}{K}
\end{equation}
where
\begin{equation}
K=ac-b^2
\end{equation}
Then easily $\omega_{12}=\widetilde{\omega}_{12}$.\\
\\
\textbf{Deffinition 3.}\\
Let $\Theta(A,B)$ be the operator 
\begin{equation}
\Theta(A,B):=\nabla_1B-\nabla_2A+q_1A+q_2B
\end{equation} 
Moreover we define the differential operator
\begin{equation}
D_{\omega}f:=\Theta(a_1f,a_2f)\textrm{, where }\omega=a_1\omega_1+a_2\omega_2.
\end{equation}
Note that $\Theta$ and $D_{w}$ can defined for vectors also.\\
\\

With the above notations notations we get the following:\\
\\
For arbitrary function $A$, we have 
\begin{equation}
\Theta\left(\nabla_1A,\nabla_2A\right)=0
\end{equation}
and
\begin{equation}
\Theta\left(-\nabla_2A,\nabla_1A\right)=\Delta_2A
\end{equation}
Also the Mainardi-Godazzi (see relations (60),(61) below) equations become
\begin{equation}
\Theta(a,b)=-\widetilde{q}_1K
\end{equation}
and
\begin{equation}
\Theta(b,c)=-\widetilde{q}_2K
\end{equation}
An interesting vector relation is
$$
\Theta\left(\overline{\epsilon}_1,\overline{\epsilon}_2\right)=0,
$$
$$
\Theta(q_1,q_2)=-K.
$$
\\
\textbf{Theorem 5.}\\
If $\omega=a_1\omega_1+a_2\omega_2$, then 
$$
d\left(f\omega\right)=\Theta(a_1f,a_2f)\omega_1\wedge\omega_2=\left(D_{\omega}f\right)\omega_1\wedge\omega_2=
$$
\begin{equation}
=\left[a_2\nabla_1f-a_1\nabla_2f+\Theta(a_1,a_2)f\right]\omega_1\wedge\omega_2
\end{equation}
and
\begin{equation}
d\omega=\Theta(a_1,a_2)\omega_1\wedge\omega_2,
\end{equation}
\begin{equation}
D_{\omega}f=
\left|\begin{array}{cc}
	\nabla_1f\textrm{ }\nabla_2f\\
	a_1\textrm{ }\textrm{ }a_2
\end{array}\right|
+\Theta(a_1,a_2)f
\end{equation}
\\
\textbf{Proof.}\\
The result follows after expanding $d(f\omega)$ and using (30).\\
\\
\textbf{Theorem 6.}\\
The diferential operator $D_{\omega}(.)$ is linear
\begin{equation}
D_{\omega}(f+g)=D_{\omega}f+D_{\omega}g
\end{equation}
The differential operator $D_{\omega}$ satisfies  
\begin{equation}
D_{\omega}(fg)=gD_{\omega}f+fD_{\omega}g-\Theta(a_1,a_2)fg
\end{equation}
\\
\textbf{Proof.}\\
The proofs of (41) and (42) are easily follow from Theorem 5.\\
\\
\textbf{Remark.}\\
The property described in relation (42) give us for $\omega_{12}=q_1\omega_1+q_2\omega_2$:  
\begin{equation}
D_{12}:=D_{\omega_{12}}=\Pi_2
\end{equation}
and
\begin{equation}
D_{12}(fg)=gD_{12}f+fD_{12}g+Kfg,
\end{equation}
where $K$ is the Gauss-curvature (from $\Theta(q_1,q_2)=-K$).\\ 
Also for all $h$ such that $\omega_0=dh=(\nabla_1h)\omega_1+(\nabla_2h)\omega_2$, we have from (34) and Theorem 6:
\begin{equation}
D_0(fg)=gD_0f+fD_0g,
\end{equation}
where $D_0:=D_{\omega_{0}}$. Also for two vectors $\overline{a},\overline{b}$ holds
$$
\Pi_2\left\langle \overline{a},\overline{b}\right\rangle=\left\langle \Pi_2(\overline{a}),\overline{b}\right\rangle+\left\langle \overline{a}, \Pi_2(\overline{b})\right\rangle+K\left\langle \overline{a},\overline{b}\right\rangle.\eqno{(45.1)}
$$
\\
\textbf{Theorem 7.}\\
We define the differential $\Theta^{III}$ as
\begin{equation}
d\widetilde{\omega}=\Theta^{III}(A_1,A_2)\omega_{31}\wedge\omega_{32}\textrm{, where }\widetilde{\omega}=A_1\omega_{31}+A_2\omega_{32}.
\end{equation}
Hence
\begin{equation}
\Theta^{III}(A_1,A_2)=\widetilde{\nabla}_1A_2-\widetilde{\nabla}_2A_1+\widetilde{q}_1A_1+\widetilde{q}_2A_2.
\end{equation}
Also we define $D^{III}_{\widetilde{\omega}}$ such
\begin{equation}
d\left(f\widetilde{\omega}\right)=\left(D^{III}_{\widetilde{\omega}}f\right)\omega_{31}\wedge\omega_{32}.
\end{equation}
Hence
\begin{equation}
D^{III}_{\widetilde{\omega}}f=\Theta^{III}(A_1f,A_2f)\textrm{, where }\widetilde{\omega}=A_1\omega_{31}+A_2\omega_{32}.
\end{equation}
\\
\textbf{Theorem 8.}
\begin{equation}
D^{III}_{\widetilde{\omega}}f=
\left|\begin{array}{cc}
	\widetilde{\nabla}_1f\textrm{  }\widetilde{\nabla}_2f\\
	A_1\textrm{  }\textrm{ }A_2
\end{array}\right|
+\Theta^{III}(A_1,A_2)f,
\end{equation}
\begin{equation}
\widetilde{\nabla}_1f=-\frac{1}{K}\left|
\begin{array}{cc}
	\nabla_1f\textrm{ }\nabla_2f\\
   b\textrm{ }\textrm{ }c 
\end{array}\right|\textrm{, }\widetilde{\nabla}_2f=\frac{1}{K}\left|
\begin{array}{cc}
\nabla_1f\textrm{ }\nabla_2f\\
a\textrm{ }\textrm{ }b
\end{array}
\right|,
\end{equation}
\begin{equation}
\Theta^{III}(A_1,A_2)=\frac{1}{K}\left(D_{\omega_{31}}A_1+D_{\omega_{32}}A_2\right),
\end{equation}
\begin{equation}
D^{III}_{\widetilde{\omega}}f=\frac{A_1}{K}
\left|\begin{array}{cc}
\nabla_1f\textrm{ }\nabla_2f\\
-a\textrm{ }\textrm{ }-b
\end{array}\right|
+\frac{A_2}{K}
\left|\begin{array}{cc}
\nabla_1f\textrm{ }\nabla_2f\\
-b\textrm{ }\textrm{ }-c
\end{array}\right|
+\Theta^{III}(A_1,A_2)f,
\end{equation}
\begin{equation}
D_{\omega}f=K\cdot D^{III}_{\widetilde{\omega}}f,
\end{equation}
and
\begin{equation}
\Theta(a_1,a_2)=K\cdot \Theta^{III}(A_1,A_2).
\end{equation}
\\
\textbf{Proof.}\\
Relation (50) is easily obtained from (49) and (47).\\
For relations (51) we have 
$$
df=\widetilde{\nabla}_1f\omega_{31}+\widetilde{\nabla}_2f\omega_{32}=\widetilde{\nabla}_1f(-a\omega_1-b\omega_2)+\widetilde{\nabla}_2f(-b\omega_1-c\omega_2)=
$$
$$
=(-a\widetilde{\nabla}_1f-b\widetilde{\nabla}_2f)\omega_1+(-b\widetilde{\nabla}_1f-c\widetilde{\nabla}_2f)\omega_2.
$$
This must be equal to $df=\nabla_1f\omega_1+\nabla_2f\omega_2$. Equating these two relations and solving with respect to $\widetilde{\nabla}_1f$ and $\widetilde{\nabla}_2f$ we get (51).\\
Relation (52) follows using the definition of $D_{\omega}f$, $\Theta^{III}(A_1,A_2)$ and the Mainardi-Godazzi equations. From (47) and (51) we get
$$
\Theta^{III}(A_1,A_2)=\widetilde{\nabla}_1A_2-\widetilde{\nabla}_2A_1+\widetilde{q}_1A_1+\widetilde{q}_2A_2=
$$ 
$$
=-\frac{1}{K}\left|
\begin{array}{cc}
	\nabla_1A_2\textrm{ }\nabla_2A_2\\
	b\textrm{ }\textrm{ }c
\end{array}
\right|-\frac{1}{K}\left|
\begin{array}{cc}
\nabla_1A_1\textrm{ }\nabla_2A_1\\
a\textrm{ }\textrm{ }b	
\end{array}
\right|+\widetilde{q}_1A_1+\widetilde{q}_2A_2=
$$
$$
=-\frac{1}{K}\left|
\begin{array}{cc}
	\nabla_1A_2\textrm{ }\nabla_2A_2\\
	b\textrm{ }\textrm{ }c
\end{array}
\right|-\frac{\Theta(b,c)A_2}{K}-\frac{1}{K}\left|
\begin{array}{cc}
\nabla_1A_1\textrm{ }\nabla_2A_1\\
a\textrm{ }\textrm{ }b	
\end{array}
\right|
-\frac{\Theta(a,b)A_1}{K}=
$$
$$
=K^{-1}\left(D_{\omega_{31}}A_1+D_{\omega_{32}}A_2\right)
$$
Hence we arrive to the result using Mainardi-Godazzi equations (60) and (61),(36),(37).\\
\\
\textbf{Theorem 9.}\\
For $\omega=a_1\omega_1+a_2\omega_2$, we set $\Theta=\Theta(a_1,a_2)$, then
\begin{equation}
\int_{\partial D}\frac{\omega}{\Theta}-\int\int_{D}
\left|\begin{array}{cc}
	\nabla_1\left(\frac{1}{\Theta}\right)  \nabla_2\left(\frac{1}{\Theta}\right)\\
	a_1 \textrm{ } \textrm{ } a_2
\end{array}\right|\omega_1\wedge\omega_2=A_{D}
\end{equation}
\\
\textbf{Corollary 3.}
\begin{equation}
\int_{\partial D}\frac{\omega_1+\omega_2}{q_1+q_2}-\int\int_{D}\left[\nabla_1\left(\frac{1}{q_1+q_2}\right)-\nabla_2\left(\frac{1}{q_1+q_2}\right)\right]\omega_1\wedge\omega_2=A_{D}
\end{equation}
\\
\textbf{Theorem 10.}
\begin{equation}
\int_{\partial D}f\omega_{31}=\int\int_{D}\left(-b\nabla_1f+a\nabla_2f+\widetilde{q}_{1}Kf\right)\omega_1\wedge\omega_2
\end{equation}
and
\begin{equation}
\int_{\partial D}f\omega_{32}=\int\int_{D}\left(-c\nabla_1f+b\nabla_2f+\widetilde{q}_2Kf\right)\omega_1\wedge\omega_2
\end{equation}
\\
\textbf{Proof.}\\The Mainardi-Godazzi equations read as
\begin{equation}
\nabla_1b-\nabla_2a+2q_2b+q_1(a-c)=0
\end{equation}
and
\begin{equation}
\nabla_1c-\nabla_2b+2q_1b-q_2(a-c)=0
\end{equation}
From relations (28) and the definition of $\Theta$ operator the result follow easily.\\
\\
\textbf{Note.}\\
If 
\begin{equation}
df=\omega_{32}\widetilde{q}_1f
\end{equation}
then
\begin{equation}
\nabla_1f=-bf\widetilde{q}_1\textrm{ and }\nabla_2f=-cf\widetilde{q}_1
\end{equation}
Hence
$$
\int_{\partial D}f\omega_{31}=\int\int_{D}(bbf\widetilde{q}_1-acf\widetilde{q}_1+\widetilde{q}_1Kf)\omega_1\wedge\omega_2=
$$
$$
=\int\int_{D}\left(-f\widetilde{q}_1K+\widetilde{q}_1Kf\right)\omega_1\wedge\omega_2=0
$$
and we get the next\\ 
\\
\textbf{Corollary 4.}\\
If $df=\omega_{32}\widetilde{q}_1f$, then
\begin{equation}
\int_{\partial D}f\omega_{31}=0
\end{equation}
If $df=-\omega_{31}\widetilde{q}_2f$, then
\begin{equation}
\int_{\partial D}f\omega_{32}=0
\end{equation}
\\

Also it holds the following formula (Liouville) 
\begin{equation}
\frac{1}{\rho_g}=\frac{d\phi}{ds}+q_1\cos\phi+q_2\sin\phi
\end{equation}
we can rewrite the above relation as
$$
\frac{1}{\rho_g}=\nabla_1\phi\frac{\omega_1}{ds}+\nabla_2\phi\frac{\omega_2}{ds}+q_1\cos\phi+q_2\sin\phi
$$
or equivalent 
\begin{equation}
\frac{1}{\rho_g}=\left(\nabla_1\phi+q_1\right)\cos\phi+\left(\nabla_2\phi+q_2\right)\sin\phi
\end{equation} 
\\
\textbf{Theorem 11.}
\begin{equation}
\frac{1}{\rho_g}=\Theta(\cos\phi,\sin\phi)
\end{equation}
\\
\textbf{Proof.}\\
From (66) and definition of $\Theta$ (relation (32)).\\
\\

If $r,\theta$ are such that
\begin{equation}
\nabla_1\phi+q_1=r\cos\theta\textrm{ and }\nabla_2\phi+q_2=r\sin\theta
\end{equation}
then we have the next\\ 
\\
\textbf{Theorem 12.}\\
If $\overline{t}=A\overline{\epsilon}_1+B\overline{\epsilon}_2$ is the tangent of a surface curve $\Gamma$, then
\begin{equation}
\left(\frac{1}{\rho_g}\right)_{\Gamma}=(A^2+B^2)\frac{\omega_{12}}{ds}+\frac{dB}{ds}A-\frac{dA}{ds}B,
\end{equation}
\\
\textbf{Proof.}\\
There hold the relations
$$
d\overline{\epsilon}_1=(q_1\omega_1+q_2\omega_2)\overline{\epsilon}_2+(a\omega_1+b\omega_2)\overline{\epsilon}_3
$$
and
$$
d\overline{\epsilon}_2=-(q_1\omega_1+q_2\omega_2)\overline{\epsilon}_1+(b\omega_1+c\omega_2)\overline{\epsilon}_3
$$
and 
$$
\frac{d\overline{t}}{d\nu}=\frac{dA}{d\nu}\overline{\epsilon}_1+A\frac{d\overline{\epsilon}_1}{d\nu}+\frac{dB}{d\nu}\overline{\epsilon}_2+B\frac{d\overline{\epsilon}_2}{d\nu}
$$
Combining the above we get the result.\\

In the special case $\overline{t}=(r\cos\theta)\overline{\epsilon}_1+(r\sin\theta)\overline{\epsilon}_2$, then
\begin{equation}
\left(\frac{1}{\rho_g}\right)_{\Gamma_1}=r^2\left(\frac{\omega_{12}}{d\nu}+\frac{d\theta}{d\nu}\right)
\end{equation}
Hence $r=1$ iff $\nu=s$ is the canonical parameter of $\Gamma_1$.\\ 
Therefore, if $\Gamma_0\in S$ is a surface curve which, its tangent at a point $P$, form with the curve $\omega_1=0$ of $S$  angle $\phi$, we can define the curve $\Gamma_1$ (of $S$) with angle $\theta$, such (69) hold and it is $r=1$ ($\nu=s$), then 
\begin{equation}
\nabla_1\phi+q_1=\cos\theta\textrm{ and }\nabla_2\phi+q_2=\sin\theta
\end{equation}   
Hence $\Gamma_1$ can be defined by (72) and also holds
\begin{equation}
\left(\frac{1}{\rho_g}\right)_{\Gamma_0}=\cos(\theta-\phi)
\end{equation}
If $\left(\frac{1}{\rho_g}\right)_{\Gamma_0}=0$, we get
\begin{equation}
\theta=\phi+\frac{\pi}{2}
\end{equation}
and relations (72) become
\begin{equation}
\nabla_1\phi+q_1=-\sin\phi\textrm{ and }\nabla_2\phi+q_2=\cos\phi.
\end{equation}
\\
\textbf{Corollary 5.}\\
For all curves $\Gamma_1$ of $\Gamma_0=\Gamma\in S$, we have
\begin{equation}
\Theta(r\cos\theta,r\sin\theta)=\Theta(q_1,q_2)=-K
\end{equation}
If $\nu=s$ is the canonical parameter of $\Gamma_1$, then
\begin{equation}
\left(\frac{1}{\rho_g}\right)_{\Gamma_1}=-K
\end{equation}
\\
\textbf{Proof.}\\
Relations (69) give
\begin{equation}
\nabla_2q_1=\nabla_2(r\cos\theta)-\nabla_2\nabla_1\phi\textrm{ and }\nabla_1q_2=\nabla_1(r\sin\theta)-\nabla_1\nabla_2\phi
\end{equation}
Hence from Gauss theorem, which is
\begin{equation}
K=\nabla_2q_1-\nabla_1q_2-q_1^2-q_2^2,
\end{equation}
we have
$$
K=\nabla_2(r\cos\theta)-\nabla_1(r\sin\theta)-\nabla_2\nabla_1\phi+\nabla_1\nabla_2\phi-
$$
$$
-q_1\left(r\sin\phi-\nabla_1\phi\right)-q_2\left(r\cos\phi-\nabla_2\phi\right)
$$
From the fundamental idenity (5), we get
$$
K=\nabla_2(r\cos\theta)-\nabla_1(r\sin\theta)
-rq_1\sin\theta-rq_2\cos\theta
$$
Also combining the above and Theorem 5, we get (77).\\
\\
\textbf{Note.}\\
If $\phi$ is the angle of a geodesic curve with $\overline{\epsilon}_1$, then
$$
\Theta(-\sin\phi,\cos\phi)=(\nabla_1\cos\phi)+(\nabla_2\sin\phi)-q_1\sin\phi+q_2\cos\phi=
$$
$$
-\sin\phi(\nabla_1\phi)+\cos\phi(\nabla_2\phi)-q_1\sin\phi+q_2\cos\phi=
$$
$$
=q_1\sin\phi+\sin^2\phi+\cos^2\phi-q_2\cos\phi-q_1\sin\phi+q_2\cos\phi=1.
$$
\\

Denote
\begin{equation}
\overline{\xi}=\cos\theta\overline{\epsilon}_1+\sin\theta\overline{\epsilon}_2
\end{equation}
and
\begin{equation}
\overline{\xi}^{*}=-\sin\theta\overline{\epsilon}_1+\cos\theta\overline{\epsilon}_2.
\end{equation}
Obviously $\left\langle \overline{\xi},\overline{\xi}^{*}\right\rangle=0$.\\ 
\\ 
\textbf{Theorem 13.}
\begin{equation}
\int_{\partial D}fdg=\int\int_{D}\left[\nabla_1f\nabla_2g-\nabla_2f\nabla_1g\right]\omega_1\wedge\omega_2
\end{equation}
\\
\textbf{Theorem 13.1}\\Assume that exists function $\lambda$ such that 
\begin{equation}
\left|\begin{array}{cc}
	\nabla_1\lambda\textrm{ }\nabla_2\lambda\\
	q_1\textrm{ }\textrm{ }q_2
\end{array}\right|=0,
\end{equation} 
then exists also function $\mu\neq 0$ such that
\begin{equation}
\int_{\partial D}\lambda d\mu=-\int_{\partial D}\mu d\lambda=\int\int_{D}K\omega_1\wedge\omega_2=\int_{\partial D}\delta\psi
\end{equation}
and 
\begin{equation}
D_{\omega_{12}}\mu=0
\end{equation}
The eigenvalue of $\mu^{\nu}$ under $D_{\omega_{12}}$ is $(\nu-1)K$ i.e
\begin{equation}
D_{\omega_{12}}\mu^{\nu}=(\nu-1)K\mu^{\nu}\textrm{, }\nu=1,2,3,\ldots
\end{equation}
\begin{equation}
D_{\omega_{12}}f(\mu)=\Pi_2f(\mu)=K\cdot(f(0)-f(\mu)+\mu f'(\mu))
\end{equation}	
\begin{equation}
\int_{\partial D}\mu\omega_{\omega_{12}}=0
\end{equation}
\begin{equation}
\Pi_2\mu=0
\end{equation}
\begin{equation}
\Pi_2\lambda=-K\lambda\textrm{, }\Pi_2\lambda^{\nu}=-K\lambda^{\nu}\textrm{, }\nu=2,3,\ldots
\end{equation}
For a smooth one valued function $f(x)$, with $f(0)=0$, we have
\begin{equation}
\Pi_2f(\lambda)=-Kf(\lambda)
\end{equation}
\\
\textbf{Proof.}\\
From (83) we get that exist $\mu\neq0$, with 
\begin{equation}
\frac{\nabla_1\lambda}{q_1}=\frac{\nabla_2\lambda}{q_2}=\frac{1}{\mu}
\end{equation}
Hence
\begin{equation}
q_1=\mu\nabla_1\lambda\textrm{, }q_2=\mu\nabla_2\lambda
\end{equation}
From Gauss theorem we have
$$
\nabla_1q_2-\nabla_2q_1+q_1^2+q_2^2+K=0
$$
Hence
$$
\nabla_1(\mu\nabla_2\lambda)-\nabla_2(\mu\nabla_1\lambda)+q_1\mu\nabla_1\lambda+q_2\mu\nabla_2\lambda+K=0
$$
or equivalent
$$
\nabla_1\mu\nabla_2\lambda+\mu\nabla_1\nabla_2\lambda-\nabla_2\mu\nabla_1\lambda-\mu\nabla_2\nabla_1\lambda+q_1\mu\nabla_1\lambda+q_2\mu\nabla_2\lambda+K=0
$$
or equivalent using the differentation identity (relation (5))
$$
\nabla_1\mu\nabla_2\lambda-\nabla_1\lambda\nabla_2\mu=-K
$$
Hence from Theorem (10) we get the result.\\
\\

Assume that
\begin{equation} \overline{\textrm{grad}(f)}=(\nabla_1f)\overline{\epsilon}_1+(\nabla_2f)\overline{\epsilon}_2
\end{equation}
and $\overline{t}=\overline{\xi}$, where $\overline{t}$ is the tangent vector of a curve $\partial D$ in the surface and $\phi$ is the angle between $\overline{t}$ and $\overline{\epsilon}_1$. Then\\
\\
\textbf{Theorem 14.}
\begin{equation}
\int_{\partial D}fd\phi=-\int\int_{D}\left[r\left\langle \overline{\textrm{grad}(f)},\overline{\xi}^{*} \right\rangle+\Pi_2f+Kf\right]\omega_1\wedge\omega_2 
\end{equation}
and
\begin{equation}
\int_{\partial D}\frac{f}{\rho_g}ds=-\int\int_{D}\left[r\left\langle \overline{\textrm{grad}(f)},\overline{\xi}^{*}\right\rangle+Kf\right]\omega_1\wedge\omega_2
\end{equation}
\\
\textbf{Proof.}\\From Theorem 13 we get
$$
d(fd\phi)=\left(\nabla_1f\nabla_2\phi-\nabla_2f\nabla_1\phi\right)\omega_1\wedge\omega_2=
$$
$$
=\left[\nabla_1f(r\sin\phi-q_2)-\nabla_2f(r\cos\phi-q_1)\right]\omega_1\wedge\omega_2=
$$ 
$$
=\left[r\nabla_1f\sin\phi-q_2\nabla_1f-r\nabla_2f\cos\phi+q_1\nabla_2f\right]\omega_1\wedge\omega_2=
$$
$$
=\left[-r\left\langle \overline{\textrm{grad}(f)},\overline{\xi}^{*}\right\rangle-(\Pi_2f+Kf)\right]\omega_1\wedge\omega_2.
$$
Also from (16) we get (95).\\
\\
\textbf{Corollary 6.}\\
If $\partial D$ is any closed geodesic curve, then
\begin{equation}
A_{D}=-\int\int_{D}r\left\langle \overline{\textrm{grad}\left(\frac{1}{K}\right)},\overline{\xi}^{*}\right\rangle\omega_1\wedge\omega_2.
\end{equation}
\\
\textbf{Lemma 3.}\\
If $\overline{t}=A\overline{\epsilon}_1+B\overline{\epsilon}_2$ is an arbitrary vector of the tangent plane of $S$ and $\omega=a_1\omega_1+a_2\omega_2$, then
$$
\frac{d\left(\overline{t}\omega\right)}{\omega_1\wedge\omega_2}=\left(D_{\omega}A-B\left|
\begin{array}{cc}
	q_1\textrm{ }q_2\\
	a_1\textrm{ }a_2
\end{array}\right|\right)\overline{\epsilon}_1+\left(D_{\omega}B+A\left|
\begin{array}{cc}
	q_1\textrm{ }q_2\\
	a_1\textrm{ }a_2
\end{array}\right|\right)\overline{\epsilon}_2+
$$
\begin{equation}
+\left(A\left|
\begin{array}{cc}
	a\textrm{ }\textrm{ }b\\
	a_1\textrm{ }a_2
\end{array}\right|+B\left|
\begin{array}{cc}
	b\textrm{ }\textrm{ }c\\
	a_1\textrm{ }a_2
\end{array}\right|\right)\overline{\epsilon}_3.
\end{equation} 
\\
\textbf{Proof.}\\
Using the identities
\begin{equation}
d\overline{\epsilon}_1=(q_1\omega_1+q_2\omega_2)\overline{\epsilon}_2+(a\omega_1+b\omega_2)\overline{\epsilon}_3
\end{equation}
\begin{equation}
d\overline{\epsilon}_2=-(q_1\omega_1+q_2\omega_2)\overline{\epsilon}_1+(b\omega_1+c\omega_2)\overline{\epsilon}_3
\end{equation}
and definition of $T_{\omega}$ we get the result after  algebraic calculations.\\
\\
\textbf{Corollary 7.}\\
If $\overline{t}=A\overline{\epsilon}_1+B\overline{\epsilon}_2$ and $\overline{t}^{*}=-B\overline{\epsilon}_1+A\overline{\epsilon}_2$, then
\begin{equation}
\left\langle d(\overline{t}\omega),d(\overline{t}^*\omega)\right\rangle=\left\langle d(\overline{t}\omega),\overline{\epsilon}_3\right\rangle \left\langle d(\overline{t}^{*}\omega), \overline{\epsilon}_3\right\rangle
\end{equation}
\textbf{Proof.}\\
Easy we get
\begin{equation}
\left\langle d\left(\overline{t}\omega\right)\times \overline{\epsilon}_3,d\left(\overline{t}^{*}\omega\right)\times \overline{\epsilon}_3\right\rangle=0
\end{equation}
From Lagrange formula
\begin{equation}
\left\langle \overline{a}\times \overline{b},\overline{c}\times \overline{d}\right\rangle=\left\langle \overline{a},\overline{c}\right\rangle \left\langle \overline{b},\overline{d}\right\rangle-\left\langle \overline{a},\overline{d}\right\rangle \left\langle \overline{b},\overline{c}\right\rangle 
\end{equation}
we get the result.\\
\\

$$
\frac{d\left(\overline{t}\omega\right)}{\omega_1\wedge\omega_2}=\frac{d\left(\overline{t}^{*}\omega\right)}{\omega_1\wedge\omega_2}\times \overline{\epsilon}_3
$$
Hence
\begin{equation}
d\left(\overline{t}\omega\right)\times\overline{\epsilon}_3=d\left(\overline{t}^{*}\omega\right)\times \overline{\epsilon}_3
\end{equation}
Therefore
$$
d\left[\left(\overline{t}-\overline{t}^{*}\right)\omega\right]\times \overline{\epsilon}_3=0
$$
Hence exists J such that 
$$
d\left[\left(\overline{t}-\overline{t}^{*}\right)\omega\right]=J\overline{\epsilon}_3.
$$
\\

Assume now that $\overline{h}=\overline{t}+C\overline{\epsilon}_3=A\overline{\epsilon}_1+B\overline{\epsilon}_2+C\overline{\epsilon}_3$, then easily as in Lemma we get\\
\\
\textbf{Theorem 15.}\\
If $\overline{h}=A\overline{\epsilon}_1+B\overline{\epsilon}_2+C\overline{\epsilon}_3$ arbitrary vector and $\omega=a_1\omega_1+a_2\omega_2$, then 
$$
\frac{d\left(\overline{h}\omega\right)}{\omega_1\wedge\omega_2}
=\left(D_{\omega}A-B\left|
\begin{array}{cc}
	q_1\textrm{ }q_2\\
	a_1\textrm{ }a_2
\end{array}\right|-C\left|
\begin{array}{cc}
	a\textrm{ }\textrm{ }b\\
	a_1\textrm{ }a_2
\end{array}\right|\right)\overline{\epsilon}_1+
$$
$$
+\left(D_{\omega}B+A\left|
\begin{array}{cc}
	q_1\textrm{ }q_2\\
	a_1\textrm{ }a_2
\end{array}\right|-C\left|
\begin{array}{cc}
	b\textrm{ }\textrm{ }c\\
	a_1\textrm{ }a_2
\end{array}\right|\right)\overline{\epsilon}_2+
$$
\begin{equation}
+\left(D_{\omega}C+A\left|
\begin{array}{cc}
	a\textrm{ }\textbf{ }b\\
	a_1\textrm{ }a_2
\end{array}\right|+B\left|
\begin{array}{cc}
	b\textrm{ }\textbf{ }c\\
	a_1\textrm{ }a_2
\end{array}\right|\right)\overline{\epsilon}_3.
\end{equation}
\\
\textbf{Proof.}\\
Use 
\begin{equation}
d\overline{\epsilon}_3=-(a\omega_1+b\omega_2)\overline{\epsilon}_1-(b\omega_1+c\omega_2)\overline{\epsilon}_2,
\end{equation}
Theorem 12 and the fact that $d(\overline{h}\omega)=d(\overline{t}\omega)+d\left(C\overline{\epsilon}_3\omega\right)$.\\
\\
\textbf{Proposition 1.}\\
If $a_1=q_1$ and $a_2=q_2$, then
\begin{equation}
\frac{d\left(\overline{h}\omega_{12}\right)}{\omega_1\wedge\omega_2}=\left(D_{12}A+C\widetilde{q}_2K\right)\overline{\epsilon}_1+\left(D_{12}B-C\widetilde{q}_1K\right)\overline{\epsilon}_2+\left(D_{12}C-K\left|\begin{array}{cc} A\textrm{ }B\\
\widetilde{q}_1\textrm{ }\widetilde{q}_2
\end{array}\right|\right)\overline{\epsilon}_3.
\end{equation}
\\
\textbf{Corollary 8.}\\
If $\omega=\omega_{12}=q_1\omega_1+q_2\omega_2$ and $\overline{t}=A\overline{\epsilon}_1+B\overline{\epsilon}_2$, then
\begin{equation}
d(\overline{t}\omega_{12})=\left((D_{12}A)\overline{\epsilon}_1+(D_{12}B)\overline{\epsilon}_2-K\left|\begin{array}{cc}
A\textrm{ }B\\
\widetilde{q}_1\textrm{ }\widetilde{q}_2
\end{array}\right|\overline{\epsilon}_3\right)\omega_1\wedge\omega_2,
\end{equation}
\begin{equation}
\Pi_2\left(\overline{t}\right)=\Pi_2(A)\overline{\epsilon}_1+(\Pi_{2}B)\overline{\epsilon}_2-K\left|\begin{array}{cc}
A\textrm{ }B\\
\widetilde{q}_1\textrm{ }\widetilde{q}_2
\end{array}\right|\overline{\epsilon}_3.
\end{equation}
\\
\textbf{Proposition 2.}\\
If $\omega=a_1\omega_1+a_2\omega_2$, then
\begin{equation}
\left\langle\frac{d\left(X\overline{\epsilon}_i\omega\right)}{\omega_1\wedge\omega_2},\overline{\epsilon}_i\right\rangle=D_{\omega}(X)\textrm{, for every } X.  
\end{equation}
Also if we set
\begin{equation}
\Omega_{ij}(X):=\left\langle\frac{d\left(X\overline{\epsilon}_i\omega\right)}{\omega_1\wedge\omega_2},\overline{\epsilon}_j\right\rangle,
\end{equation} 
then 
\begin{equation}
\Omega_{ii}(X)=D_{\omega}(X),
\end{equation}
\begin{equation}
\Omega_{ij}\left(X\right)=X\frac{\left\langle \left(d\overline{\epsilon}_i\wedge \omega\right),\overline{\epsilon}_j\right\rangle}{\omega_1\wedge \omega_2}\textrm{, when }i\neq j.
\end{equation}
\\

Also if $\overline{h}=X_1\overline{\epsilon}_1+X_2\overline{\epsilon}_2+X_3\overline{\epsilon}_3$, then
$$
\left\langle\frac{d\left(\overline{h}\omega\right)}{\omega_1\wedge\omega_2},\overline{\epsilon}_j\right\rangle=\sum^{3}_{i=1}\Omega_{ij}(X_i)
$$
Hence
$$
\left\langle\frac{d\left(\overline{h}\omega\right)}{\omega_1\wedge\omega_2},\overline{\mu}\right\rangle=\sum^{3}_{j=1}\Omega_{ij}\left(X_i\right)\mu_j.
$$
Now define
$$
D_{\omega}\left(\overline{h}\right):=\sum^{3}_{j=1}D_{\omega}(h_j)\overline{\epsilon}_{j}.
$$
Then 
$$
d\left\langle \overline{h}\omega,\overline{\mu}\right\rangle=\sum^{3}_{i=1}
\left|\begin{array}{cc}
\nabla_1\left(\mu_ih_i\right)\textrm{  }\nabla_2\left(\mu_i h_i\right)\\
a_1\textrm{  }a_2\\
\end{array}\right|+\Theta(a_1,a_2)\sum^{3}_{i=1}h_i\mu_i=
$$
$$
=\sum^{3}_{i=1}\left(h_i\left|\begin{array}{cc}
\nabla_1\mu_i\textrm{  }\nabla_2 \mu_i\\
a_1\textrm{  }a_2
\end{array}\right|+\mu_i\left|\begin{array}{cc}
\nabla_1 h_i\textrm{  }\nabla_2 h_i\\
a_1\textrm{  }a_2
\end{array}\right|+\Theta(a_1,a_2)\mu_ih_i\right)=
$$
$$
=\left\langle D_{\omega}\overline{h},\overline{\mu}\right\rangle+\left\langle \overline{h},D_{\omega}\overline{\mu}\right\rangle-\Theta(a_1,a_2)\left\langle \overline{h},\overline{\mu}\right\rangle.
$$
Hence we get the next:\\
\\
\textbf{Proposition 3.}
$$
d\left\langle \overline{h}\omega,\overline{\mu}\right\rangle=\left\langle \overline{h},D_{\omega}\overline{\mu}\right\rangle+\left\langle D_{\omega}\overline{h},\overline{\mu}\right\rangle-\Theta(a_1,a_2)\left\langle \overline{h},\overline{\mu}\right\rangle.
$$
\\
\textbf{Proposition 4.}\\
If we set $\Gamma_{12}:=\left|
\begin{array}{cc}
	q_1\textrm{ }q_2\\
	a_1\textrm{ }a_2
\end{array}\right|$, $\Gamma_{13}:=\left|
\begin{array}{cc}
	a\textrm{ }\textrm{ }b\\
	a_1\textrm{ }a_2
\end{array}\right|$, $\Gamma_{23}:=\left|
\begin{array}{cc}
	b\textrm{ }\textrm{ }c\\
	a_1\textrm{ }a_2
\end{array}\right|$, then
\begin{equation}
\Omega_{12}(X)=X\Gamma_{12}\textrm{,  }\Omega_{13}(X)=X\Gamma_{13}\textrm{, }\Omega_{23}(X)=X\Gamma_{23},
\end{equation}
\begin{equation}
\Omega_{ij}(X)=-\Omega_{ji}(X)\textrm{, }i\neq j.
\end{equation}
\\
\textbf{Proof.}\\
Easy from Proposition 2.\\
\\
\textbf{Theorem 16.}\\
If $\overline{X}=X_1\overline{\epsilon}_1+X_2\overline{\epsilon_2}+X_3\overline{\epsilon}_3$, then
\begin{equation}
\frac{d\left(\overline{X}\omega\right)}{\omega_1\wedge \omega_2}=\left(D_{\omega}X_1-X_2\left|\begin{array}{cc}
	q_1\textrm{ }\textrm{ }q_2\\
	a_1\textrm{ }a_2
\end{array}\right|-X_3\left|\begin{array}{cc}
	a\textrm{ }\textrm{ }b\\
	a_1\textrm{ }a_2
\end{array}\right|\right)\overline{\epsilon}_1+
$$
$$
+\left(D_{\omega}X_2+X_1\left|\begin{array}{cc}
	q_1\textrm{ }\textrm{ }q_2\\
	a_1\textrm{ }a_2
\end{array}\right|-X_3\left|\begin{array}{cc}
	b\textrm{ }\textrm{ }c\\
	a_1\textrm{ }a_2
\end{array}\right|\right)\overline{\epsilon}_2+
$$
$$
+\left(D_{\omega}X_3+X_1\left|\begin{array}{cc}
	a\textrm{ }\textrm{ }b\\
	a_1\textrm{ }a_2
\end{array}\right|+X_2\left|\begin{array}{cc}
	b\textrm{ }\textrm{ }c\\
	a_1\textrm{ }a_2
\end{array}\right|\right)\overline{\epsilon}_3.
\end{equation}
\\
\textbf{Proof.}
$$
\frac{d\left(\overline{X}\omega\right)}{\omega_1\wedge \omega_2}=\sum^{3}_{i=1}D_{\omega}(X_i)\overline{\epsilon}_i+\sum^{3}_{i<j}\Omega_{ij}(X_i)\overline{\epsilon}_j+\sum^{3}_{i>j}\Omega_{ij}(X_i)\overline{\epsilon}_j=
$$
$$
=\sum^{3}_{i=1}D_{\omega}(X_i)\overline{\epsilon}_i+\sum^{3}_{i<j}X_i\Gamma_{ij}\overline{\epsilon}_j+\sum^{3}_{i>j}X_{i}\Gamma_{ij}\overline{\epsilon}_j=
$$
$$
=\sum^{3}_{i=1}D_{\omega}(X_i)\overline{\epsilon}_i+\sum^{3}_{i<j}X_i\Gamma_{ij}\overline{\epsilon}_j+\sum^{3}_{i<j}X_{j}\Gamma_{ji}\overline{\epsilon}_i=
$$
$$
=\sum^{3}_{i=1}D_{\omega}(X_i)\overline{\epsilon}_i+\sum^{3}_{i<j}X_i\Gamma_{ij}\overline{\epsilon}_i-\sum^{3}_{i<j}X_{j}\Gamma_{ij}\overline{\epsilon}_i.
$$
Hence using Propositions 2,4 we get the result.\\
\\
\textbf{Corollary 9.}\\
If $\overline{h}=A\overline{\epsilon}_1+B\overline{\epsilon}_2+C\overline{\epsilon}_3$ and $K=0$, then
\begin{equation}
\frac{d\left(\overline{h}\omega_{12}\right)}{\omega_1\wedge\omega_2}=\left|
\begin{array}{cc}
	\nabla_1A\textrm{  }\nabla_2A\\
	q_1\textrm{  }\textrm{  }q_2
\end{array}\right|\overline{\epsilon}_1+\left|
\begin{array}{cc}
	\nabla_1B\textrm{  }\nabla_2B\\
	q_1\textrm{  }\textrm{  }q_2
\end{array}\right|\overline{\epsilon}_2+\left|
\begin{array}{cc}
	\nabla_1C\textrm{  }\nabla_2C\\
	q_1\textrm{  }\textrm{  }q_2
\end{array}\right|\overline{\epsilon}_3.
\end{equation}
\\
\textbf{Proposition 5.}\\
Let $S$ be a surface and $\overline{t}=\frac{d\overline{x}}{ds}=\cos(\phi)\overline{\epsilon}_1+\sin(\phi)\overline{\epsilon}_2$ be the tangent vector of surface curve $\Gamma$. If $\phi$ is the angle between $\overline{t}$ and $\overline{\epsilon}_1$ and $\frac{\pi}{2}-\phi$ is the angle between $\overline{t}$ and $\overline{\epsilon}_2$, then
\begin{equation}
\frac{d\left(\overline{t}\omega_{12}\right)}{\omega_1\wedge\omega_2}=
\left|\begin{array}{cc}
\nabla_1\phi\textrm{ }\nabla_2\phi\\
q_1\textrm{ }q_2 	
\end{array}\right|
\overline{t}^{*}-K\overline{t}+\left(\cos(\phi) \left|\begin{array}{cc}
a\textrm{ }\textrm{ }b\\
q_1\textrm{ }\textrm{ }q_2 	
\end{array}\right|+\sin(\phi)\left|\begin{array}{cc}
b\textrm{ }\textrm{ }c\\
q_1\textrm{ }q_2 	
\end{array}\right|\right)\overline{\epsilon}_3.
\end{equation}
Hence
\begin{equation}
\Pi_2\left(\overline{t}\right)=
\left|\begin{array}{cc}
\nabla_1\phi\textrm{ }\nabla_2\phi\\
q_1\textrm{ }q_2 	
\end{array}\right|
\overline{t}^{*}-K\overline{t}+\left(\cos(\phi) \left|\begin{array}{cc}
a\textrm{ }\textrm{ }b\\
q_1\textrm{ }\textrm{ }q_2 	
\end{array}\right|+\sin(\phi)\left|\begin{array}{cc}
b\textrm{ }\textrm{ }c\\
q_1\textrm{ }q_2 	
\end{array}\right|\right)\overline{\epsilon}_3,
\end{equation}
where $\overline{t}^{*}=(-\sin\phi)\overline{\epsilon}_1+(\cos\phi)\overline{\epsilon}_2$.\\
\\
\textbf{Proof.}\\
Seting the parameters of the problem to the formula of Proposition 2, we easily get the result.\\
\\

\centerline{\bf References}

[1]: N.K. Stephanidis. ''Differential Geometry  I''. Zitis Publications. Thessaloniki, 1995.

[2] : N.K. Stephanidis. ''Differential Geometry II''. Zitis Publications. Thessaloniki, 1987.

[3]: Nirmala Prakash. ''Differential Geometry. An Integrated Approach''. Tata McGraw-Hill Publishing Company Limited. New Delhi. 1981.

\end{document}